\documentclass[12pt]{article}
\usepackage{e-jc}

\usepackage{amsmath}
\usepackage{amsthm}
\usepackage{enumerate}

\newtheoremstyle{theorem_style}%
    {8pt plus2pt minus4pt}%
    {8pt plus2pt minus4pt}%
    {\itshape}%
    {}%
    {\bfseries\scshape}%
    {}%
    {6pt}
    {}%

\newtheoremstyle{remark_style}%
    {8pt plus2pt minus4pt}%
    {8pt plus2pt minus4pt}%
    {\upshape}
    {}%
    {\bfseries\scshape}%
    {}%
    {6pt}
    {}%

\theoremstyle{theorem_style}
\newtheorem{theorem}{Theorem}
\newtheorem{lemma}[theorem]{Lemma}

\theoremstyle{remark_style}
\newtheorem{remark}[theorem]{Remark}

\def\E{\mathbb{E}}
\def\Var{{Var\;}}
\def\Det{\mbox{Det}}

\def\a{\alpha} \def\b{\beta}  
\def\e{\epsilon}    
\def\G{\Gamma}  \def\k{\kappa}
     
\def\La{\Lambda}  \def\n{\nu} \def\p{\pi}\def\x{\xi}
   \def\S{\Sigma}

\def\cY{{\cal Y}}

\def\cB{{\cal B}}

\def\cX{{\cal X}}

\newcommand{\brac}[1]{\left(#1\right)}

\newcommand{\bfrac}[2]{\left(\frac{#1}{#2}\right)}

\def\cE{{\cal E}}

\newcommand{\beq}[1]{\begin{equation}\label{#1}}
\newcommand{\eeq}{\end{equation}}

\newcommand{\set}[1]{\left\{#1\right\}}

\def\E{\mbox{{\bf E}}}
\def\Var{\mbox{{\bf Var}}}
\def\Pr{\mbox{{\bf Pr}}}
\def\whp{{\bf whp}}

\newcommand{\cA}{{\cal A}}

\def\cU{{\cal U}}

\def\E{\mbox{{\bf E}}}
\def\Var{\mbox{{\bf Var}}}
\def\Pr{\mbox{{\bf Pr}}}
\def\whp{{\bf whp}}

\newcommand{\Hp}[1]{H_{n,#1;k}}
\def\bX{Y}

\def\hG{\hat{\G}}

\newcommand{\scr}[2]{_{#1}^{(#2)}}

\title{Loose Hamilton Cycles in Random Uniform Hypergraphs}
\author{
Andrzej Dudek and
Alan Frieze\thanks{Supported in part by NSF grant CCF0502793.}\\
Department of Mathematical Sciences\\[-0.8ex]
Carnegie Mellon University\\[-0.8ex]
Pittsburgh, PA 15213\\[-0.8ex]
}

\date{\today}

\begin{document}
\maketitle

\begin{abstract}
In the random $k$-uniform hypergraph $H_{n,p;k}$ of order~$n$
each possible $k$-tuple appears independently with probability~$p$.
A loose Hamilton cycle is a cycle of order~$n$ in which every pair
of adjacent edges intersects in a single vertex. We prove that if
$p n^{k-1}/\log n$ tends to infinity with~$n$ then
$$
\lim_{\substack{n\to \infty\\2(k-1) |n}}\Pr(H_{n,p;k}\ contains\ a\ loose\ Hamilton\ cycle)=1.
$$
This is asymptotically best possible.
\end{abstract}

\section{Introduction}
The threshold for the existence of Hamilton cycles in the random graph $G_{n,p}$ has been known
for many years, see, e.g., \cite{AKS}, \cite{Boll} and \cite{KS}.
There have been many generalizations
of these results over the years and the problem is well understood. It is natural to try to extend
these results to hypergraphs and this has proven to be difficult. The famous P\'osa lemma fails to
provide any comfort and we must seek new tools. In the graphical case, Hamilton cycles and perfect
matchings go together and our approach will be to build on the
deep and difficult result of Johansson,
Kahn and Vu \cite{JKV}, as well as what we have learned from the graphical case.

A {\em $k$-uniform hypergraph} is a pair $(V,E)$ where $E\subseteq \binom{V}{k}$.
In the random $k$-uniform hypergraph $H_{n,p;k}$ of order~$n$ each
possible $k$-tuple appears independently with probability~$p$.
We say that a $k$-uniform hypergraph $(V,E)$ is a {\em loose Hamilton cycle} if there exists a cyclic
ordering of the vertices $V$ such that every edge consists of $k$ consecutive vertices and every
pair of consecutive edges intersects in a single vertex. In other words, a loose Hamilton cycle has
the minimum possible number of edges among all cycles on $|V|$ vertices.
In a recent paper the second author proved the following:
\begin{theorem}[Frieze \cite{F}]\label{th1a}
There exists an absolute constant $K>0$ such that if $p\geq K(\log n)/n^2$ then
$$\lim_{\substack{n\to \infty\\4 |n}}\Pr(H_{n,p;3}\ contains\ a\ loose\ Hamilton\ cycle)=1.$$
\end{theorem}

In this paper we refine the above theorem to $k\geq 4$. Here we state our main result.
\begin{theorem}\label{th1}
Let $k\ge 3$. If $p n^{k-1}/\log n$ tends to infinity together with~$n$ then
$$\lim_{\substack{n\to \infty\\2(k-1) |n}}\Pr(H_{n,p;k}\ contains\ a\ loose\ Hamilton\ cycle)=1.$$
\end{theorem}
Thus $(\log n)/ n^{k-1}$ is the asymptotic threshold for the
existence of loose Hamilton cycles, at least for $n$ a multiple of $2(k-1)$. This is because
if $p\leq (1-\e) (k-1)! (\log n) / n^{k-1}$ and $\e>0$ is
constant, then \whp\footnote{An event $\cE_n$
occurs {\em with high probability}, or \whp\ for brevity, if
$\lim_{n\rightarrow\infty}\Pr(\cE_n)=1$.} $\Hp{p}$ contains isolated vertices.

Notice that the necessary divisibility requirement for a
$k$-uniform hypergraph to have a loose Hamilton
cycle is $(k-1) | n$. In our approach we needed to
assume more, namely, $2(k-1) | n$ (the same is true for Theorem~\ref{th1a}).

There are other ways of defining Hamilton cycles in
hypergraphs, depending on the size of the intersection
of successive edges. As far as we know, when these intersections have more than one vertex, nothing
significant is known about existence thresholds.

Our proof uses a second moment calculation on a related
problem. We cannot apply a second moment calculation directly to the number
of Hamilton cycles in $H_{n,p;k}$, this does not work.

\section{Proof of Theorem \ref{th1}}

Fix an integer~$k\ge 3$. Set $\k=k-2$ and let $n=2(k-1)m$.
We immediately see the divisibility requirement $2(k-1)|n$.
Let $p n^{k-1}/\log n$ tend to infinity together with~$n$ (or
equivalently together with~$m$). From on now, all asymptotic notations are with respect to $m$.

We start with a special case of the theorem of \cite{JKV}. Let $S$ and $T$ be disjoint sets.
Let $\G=\G(S,T,p)$ be the random $k$-uniform hypergraph
such that each $k$-edge in $\binom{S}{2}\times \binom{T}{\k}$ is independently
included with probability $p$.
Assuming that $|S|=2m$ and $|T|=\k m$ for some positive integer $m$,
a \textit{perfect matching} of $\G$ is a set of $m$ $k$-edges
$\{s_{2i-1},s_{2i},t_{i,1},\ldots,t_{i,\k}\}$, $1\le i\le m$,
such that $\set{s_1,\ldots,s_{2m}}=S$ and $\set{t_{1,1},\ldots,t_{m,\k}}=T$.
\begin{theorem}[Johansson, Kahn and Vu~\cite{JKV}]\label{th2}
There exists an absolute constant $K>0$ such that if $p\ge K(\log n)/n^{k-1}$ then
\whp\ $\G$ contains a perfect matching.
\end{theorem}
This version is not actually proved in \cite{JKV},
but can be obtained by straightforward changes to their proof.

Now we (deterministically) partition $[n]$ into $X=[2m]$ and
$\bX=[2m+1,n]$, where clearly $|X|=2m$ and $|\bX|=2\k m$.
We show that $\G(X,\bX,p)$, which can be viewed as the subgraph of $H_{n,p;k}$ induced by
$\binom{X}{2}\times \binom{\bX}{\k}$, contains a loose Hamilton cycle \whp.
Such a Hamilton cycle will consist of $2m$ edges of the form
$\{x_{i},x_{i+1},y_{i,1},\ldots,y_{i,\k}\}$, where $1\le i\le 2m$, $x_{2m+1}=x_1$,
$\set{x_1,\ldots,x_{2m}}=X$ and $\set{y_{1,1},\ldots,y_{2m,\k}}=\bX$.

Let $d$ be an arbitrarily large even positive integer constant.
Let $\cX$ be a set of size $2dm$ representing
$d$ copies of each $x\in X$. Denote the $j$th copy of $x\in X$ by $x\scr{}{j}\in \cX$
and let $\cX_x=\set{x\scr{}{j},\,j=1,2,\ldots,d}$.
Then let $X_1,X_2,\ldots,X_{d}$ be a uniform random partition of $\cX$ into $d$ sets of size $2m$.
Define $\psi_1:\cX\to X$ by $\psi_1(x\scr{}{j})=x$ for all $j$ and $x\in X$.
Similarly, we let $\cY$ be a set of size $d\k m$ representing
$d/2$ copies of each $y\in \bX$. Denote the $j$th copy of $y\in \bX$ by $y\scr{}{j} \in \cY$
and let $\cY_y=\set{y\scr{}{j},\,j=1,2,\ldots,d/2}$.
Then let $Y_1,Y_2,\ldots,Y_{d}$ be a uniform random partition of $\cY$ into $d$ sets of size $\k m$.
Define $\psi_2:\cY\to \bX$ by $\psi_2(y\scr{}{j})=y$ for all $y\in\bX$.
Finally, let $\psi:\binom{\cX}{2}\times \binom{\cY}{\k} \to X^2 \times \bX^{\k}$
be such that $\psi(\n_1,\n_2,\x_1,\x_2,\dots,\x_\k) =
(\psi_1(\n_1),\psi_1(\n_2),\psi_2(\x_1),\psi_2(\x_2),\dots,\psi_2(\x_\k))$.

Define $p_1$ by $p=1-(1-p_1)^\a$ where $\a=e^{2\k d}$.
With this choice, we can generate $H_{n,p;k}$ as the union of $\a$ independent copies of
$H_{n,p_1;k}$. Similarly, define $p_2$ by $p_1=1-(1-p_2)^d$. Finally define
$p_3$ by $p_2=1-(1-p_3)^\b$ where $\b=d^2(d/2)^\kappa$.
Observe that $p_in^{k-1}/\log n\to\infty$ for $i=1,2,3$ as $n\to\infty$.
In this way, $H_{n,p;k}$ is represented as the union of $d\a\b$ independent copies of $H_{n,p_3;k}$.

Now let an edge $\set{\n_1,\n_2, \x_1,\x_2,\ldots,\x_\k}$ of  $\G(X_j,Y_j,p_2)$, 
$1\le j\le d$, be {\em spoiled}
if $\psi_1(\n_1)=\psi_1(\n_2)$ or
there exist $1\leq r<s\leq \k$ such that $\psi_2(\x_r)=\psi_2(\x_s)$.
Let $\hG(X_j,Y_j,p_2)$ be obtained from $\G(X_j,Y_j,p_2)$ by removing all spoiled edges.

As we already mentioned $H_{n,p;k}$ is represented as the union
of $d\a\b$ independent copies of $H_{n,p_3;k}$.
We group the $d\a\b$ copies of $H_{n,p_3;k}$ together into $\a$
sets $\cA_1,\cA_2,\ldots,\cA_\a$ in such a way that each collection
$\cA_i$, $1\le i\le \a$, consists of $d$ sub-collections $\cB_{i,j}$,
$1\le j\le d$, where $\cB_{i,j}$ comprises
$\b$ independent copies of $H_{n,p_3;k}$.
Let $\La_{i,j}$ denote the union of these $\b$ copies in
$\cB_{i,j}$ and let $\S_i$ denote the union of $\La_{i,j}$ over all $1\le j\le d$.
Basically $\La_{i,j}$ and $\S_i$ can be viewed as copies
of $H_{n,p_2;k}$ and $H_{n,p_1;k}$, respectively.

Now for fixed  $1\le i\le \a$ and $1\le j\le d$, we couple an
independent copy of $\hG(X_j,Y_j,p_2)$ with
a sub-hypergraph (induced by $\binom{X}{2}\times \binom{\bX}{\k}$)
of the union of $\b$ independent copies of $H_{n,p_3;k}$ in $\cB_{i,j}$
as follows. First we enumerate these $\b$ copies of $H_{n,p_3;k}$ as $H_{j_1,\ldots,j_k}$,
where $1\le j_1, j_2\le d$ and $1\le j_3,\ldots,j_k\le d/2$.
Next we place $\{x_1<x_2,y_1<y_2<\cdots<y_\k\}$ in $H_{j_1,\ldots,j_k}$,
whenever there exist $j_1,\ldots,j_k$ such that
$\{x\scr{1}{j_1},x\scr{2}{j_2},y\scr{1}{j_3},\ldots,y\scr{\k}{j_k}\}$ is
an edge in $\hG(X_j,Y_j,p_2)$.

Fix $1\le i\le \a$ for the moment and consider $\La_{i,j}$ for all $1\le j \le d$.
Let $M_{j}$, $1\le j\le d$, be a perfect matching of $\G(X_j,Y_j,p_2)$
as promised by Theorem \ref{th2}.
At this point what we can say is that $X_1,X_2,\ldots,X_d$
is a uniform random partition of $\cX$ and $Y_1,Y_2,\ldots,Y_d$
is a uniform random partition of $\cY$. Furthermore, if $M_j$ exists then by symmetry we
can assume that it is a uniformly random matching of $\G(X_j,Y_j,p_2)$. What we want though
are unspoiled matchings.
Fortunately, it is reasonably likely that $M_j$ contains no spoiled edges.
Our argument will be (see Lemma \ref{spoil} below) that there is a probability of at least
$e^{-\k d}$ that $M_{j}\subseteq \hG(X_j,Y_j,p_2)$ simultaneously
for all $1\le j\le d$. That means that with the same probability
$\psi(M_{j}) \subseteq \La_{i,j}$ simultaneously for all $1\le j\le d$,
i.e., $\psi(M_1\cup M_2\cup\cdots\cup M_d)\subseteq \S_i$.
It follows that
then with probability at least
\beq{tue}
1-((1-o(1))(1-e^{-\k d}))^\a \ge 1-e^{-e^{\k d}}
\eeq
there is an $i$ such that $\S_i$
contains a copy of the following hypergraph $\La_d=\psi(M_1\cup M_2\cup\cdots\cup M_d)$,
where each $M_j$ is a random perfect matching
of $\hG(X_j,Y_j,p_1)$, i.e., $M_j$ has no spoiled edges. (The first $1-o(1)$ factor in \eqref{tue}
comes from the use of Theorem \ref{th2}). We will choose such an $i$ for constructing 
$\La_d$. These matchings
are still independently chosen, once we have fixed the
partitions $X_1,X_2,\ldots,X_d$ and $Y_1,Y_2,\ldots,
Y_d$ and each $M_j$ is uniformly random from $\hG(X_j,Y_j,p_1)$ by symmetry.
On the other hand, the partitions of $\cX,\cY$ are no longer uniform. Their 
probability of selection depends on
how many unspoiled matchings they contain.

Our main auxiliary result, see Theorem \ref{thcolor},
shows that the hypergraph $\La_d$ contains a loose
Hamilton cycle with probability
at least $1-{3\k}/{d}$.
Because we have $pn^{k-1}/\log n\to\infty$ we can make $d$ arbitrarily large and
consequently this and \eqref{tue} imply that
$$\lim_{n\to\infty}\Pr(H_{n,p;k}\text{ has no Hamilton cycle})
\leq \lim_{d\to\infty}\left(e^{-e^{\k d}} + \frac{3\k}{d}\right)=0.$$
This completes the proof of Theorem~\ref{th1}.

\begin{remark}\label{rem1}
It is important to understand the distribution of $\La_d$. It is the union of matchings
$M_1,M_2,\ldots,M_d$ obtained by repeating the following experiment until the occurrence
of $\cU$:
\begin{enumerate}[(i)]
\item choose uniform random partitions of $\cX,\cY$; and then
\item choose uniform random matchings $M_j$ of $\G(X_j,Y_j,p_2)$.
\end{enumerate}
Lemma \ref{spoil} shows that we should not have to wait too long until $\cU$ occurs.
We do not choose one set of partitions and then choose the
matchings conditional on $\cU$.
\end{remark}
\section{Auxiliary results}
We will use a configuration model type of construction to analyze $\La_d$
(see, e.g., \cite{B} or
Section~9.1 in~\cite{JLR}). $\cX$ is represented as $2dm$ {\em points} partitioned into $2m$ cells
$\cX_x,x\in X$ of
$d$ points. Analogously $\cY$ is represented as $d\k m$ points partitioned into
$2\k m$ cells $\cY_y,y\in\bX$ of $d/2$ points.
To construct $\La_d$ we
take a random pairing of $\cX$ into $dm$ sets $e_1,e_2,\ldots,e_{dm}$ of size two and a random
partition $f_1,f_2,\ldots,f_{dm}$ of $\cY$
into $dm$ sets of size $\k$.
The edges of $\La_d$ will be $\psi(e_\ell\cup f_\ell)$ for $\ell=1,2,\ldots,md$.
We condition on $\cU$.

We will now argue that this model is justified. First of all ignore the event $\cU$.
To generate $M_1,M_2,\ldots,M_d$,
we can take a random
permutation $\p_1$ of $\cX$ and a random permutation $\p_2$ of $\cY$. We let
$X_j=\set{\p_1(2(j-1)m+i),\,i=1,\ldots,2m}$ and then $M_{j,X}$
will consist of $e_\ell=\set{\p_1(2\ell-1),\p_1(2\ell)}$ for $\ell=(j-1)m+1,\ldots,jm$.
We construct the $f_\ell$ and $Y_j$ and $M_{j,Y}$ in a similar way from $\p_2$.
So $\p_1,\p_2$ generate the same hypergraph when
viewed either as originally described in terms of $M_1,M_2,\ldots,M_d$ or as described in terms of a
configuration model. Each sequence $M_1,M_2,\ldots,M_d$ is equally likely in both models.
The relationship between models will therefore
continue to hold even if we condition on the event $\cU$.

As already noted in Remark \ref{rem1}, $\La_d$ is the above model conditioned on the event
$\cU$. We generate a conditioned sample by repeatedly generating $M_1,M_2,\ldots,M_d$ until
the event $\cU$ occurs.
In our analysis of the configuration model we deal with $\cU$ directly. We use a second
moment method and compute our moments conditional on $\cU$.

\subsection{Spoiled edges}

Suppose that for every $1\le j\le d$ there exists a perfect matching
$M_{j}$ of $\G(X_j,Y_j,p_2)$. We show that it is reasonably likely
that $M_1\cup\dots\cup M_d$ contains no spoiled edges.

Let $\cU$ be the event:
$$\{M_j\subseteq \hG(X_j,Y_j,p_2),\text{for each }j=1,2,\ldots,d\}=
\{M_1\cup\dots\cup M_d\text{ contains no spoiled edges}\}.$$

\begin{lemma}\label{spoil}
Suppose that ${\k}\ge 1$ and $d$ is a positive even integer. Then,
\footnote{We write $A_m\sim B_m$
to signify that $A_m=(1+o(1))B_m$ as $m\to\infty$.}
$$
\Pr(\cU\mid M_j
\text{ exists for each } j=1,2,\ldots,d)
\sim\exp\set{-\frac{d-1}{2}-\frac{(\k-1)(d-2)}{4}} \geq e^{-\k d}.
$$
\end{lemma}
\begin{proof}
Our model for $M_j$ will be a collection of sets
$\set{x_{j,2\ell-1},x_{j,2\ell},Z_{j,\ell}}$,
where $M_{j,X}=$\\ $\set{x_{j,1},x_{j,2}},\ldots,\set{x_{j,2m-1}x_{j,2m}}$ is a random pairing
of $X_j$ and $M_{j,\bX}=Z_{j,1},Z_{j,2},\ldots,Z_{j,m}$ is a random partition of $Y_j$ into
sets of size~$\k$.
We can obtain all of the $\set{x_{j,2\ell-1},x_{j,2\ell}}$, for all $j$ and $\ell$,
by taking a random permutation of $\cX$ and then considering it in $dm$ consecutive
sub-sequences $I_1,I_2,\ldots,I_{dm}$ of length $2$.
Let $S_1$ denote the number of pairs $\n_1,\n_2$ of
elements in $\cX$ with $\psi_1(\n_1)=\psi_1(\n_2)$ that appear in some
$I_\ell$.
Similarly, we can obtain all of the the $Z_{j,\ell}$ by taking a random
permutation of $\cY$ and then considering it in $dm$
consecutive sub-sequences $J_1,J_2,\ldots,J_{dm}$ of length $\k$. Let
now $S_2$ denote the number of pairs $\x_1,\x_2$ of elements in $\cY$
with $\psi_2(\x_1)=\psi_2(\x_2)$ that appear in some $J_\ell$.
Then for any constant $t\geq 1$, we obtain
$$
\E(S_1(S_1-1)\cdots(S_1-t+1))\sim t!\binom{dm}{t}\brac{\frac{d-1}{2dm-O(1)}}^t
\sim\bfrac{d-1}{2}^t,
$$
and
$$
\E(S_2(S_2-1)\cdots(S_2-t+1))\sim t!\binom{dm}{t}\brac{\binom{\k}{2}\frac{d/2-1}{d\k m-O(1)}}^t
\sim\bfrac{(\k-1)(d-2)}{4}^t.
$$
It follows that $S_1$ and $S_2$ are asymptotically Poisson with means
$\frac{(d-1)}{2}$ and $\frac{(\k-1)(d-2)}{4}$, respectively. Now $S_1$ and $S_2$ are
independent and so
$S_1+S_2$ is asymptotically Poisson with mean
$\frac{(d-1)}{2} + \frac{(\k-1)(d-2)}{4}$ and
\begin{align*}
\Pr(M_j\subseteq \hG(X_j,Y_j,p_2),&\text{ for each }j=1,2,\ldots,d\mid M_j
\text{ exists for each } j=1,2,\ldots,d)\\
&= \Pr(S_1+S_2=0 \mid M_j \text{ exists for each } j=1,2,\ldots,d)\\
&\sim\exp\set{-\frac{d-1}{2}-\frac{(\k-1)(d-2)}{4}}\\
&\geq e^{-\k d},
\end{align*}
as required.
\end{proof}

\subsection{Loose Hamilton cycles in random bipartite hypergraphs}

Recall that $\cX$ is a set of size $2dm$ representing $d$ copies
of each $x\in X$ and $\cY$ is a set of size $d\k m$ representing $d/2$
copies of each $y\in \bX$, where $|X|=2m$ and $|\bX|=2\k m$. Let $X_1,X_2,\ldots,X_{d}$
be a uniform random partition of $\cX$ into $d$ sets of size $2m$ and let $Y_1,Y_2,\ldots,Y_{d}$
be a uniform random partition of $\cY$ into $d$ sets of size $\k m$. For every $1\le j\le d$,
let $M_j$ be a random matching of $\binom{X_j}{2}\times \binom{Y_j}{\k}$ conditioned on $\cU$
i.e. without spoiled edges.
That means $M_j$ is a set of $m$ disjoint $k$-edges in $\binom{X_j}{2}\times \binom{Y_j}{\k}$
such that no edge contains two representatives of the same element of $X\cup \bX$.
Let $\La_d = \psi(M_1\cup\dots\cup M_d)$.

\begin{theorem}\label{thcolor}
Suppose that ${\k}\ge 1$ and $d$ is a sufficiently large positive even integer. Then,
$$
\Pr(\La_d\text{ contains a loose Hamilton cycle})\geq 2-(1+o(1))\sqrt{\frac{d}{d-2(\k+1)}}
\geq 1-\frac{3\k}{d}.
$$
\end{theorem}
A similar result for $\k=1$ was already established by Janson and Wormald~\cite{JW}
using a different terminology.

Let $H$ be a random variable which counts the number of loose Hamilton cycles in $\La_d$
such that the edges only intersect in $X$. Note that every such loose Hamilton cycle
induces an ordinary Hamilton cycle of length~$2m$ in $X$ and a partition of $\bX$ into $\k$-sets.

\begin{lemma}\label{lem:expectation}
Suppose that ${\k}\ge 1$ and $d$ is a positive even integer. Then,
$$
\E(H) \sim e^{(\k+1)/2} \pi \sqrt{\frac{{\k}(d-2)}{d}}
\left( \frac{(d-1)(d-2)^{\frac{{\k}+1}{2}(d-2)}}{d^{\frac{{\k}+1}{2}(d-2)}} \right)^{2m}.
$$
Hence, $\lim_{m\to\infty} \E(H) = \infty$ for every $d > e^{{\k}+1}+1$.
\end{lemma}

\noindent
The last conclusion holds since for $d> e^{{\k}+1}+1$,
\begin{eqnarray*}
\frac{(d-1)(d-2)^{\frac{{\k}+1}{2}(d-2)}}{d^{\frac{{\k}+1}{2}(d-2)}}
&=&(d-1)\brac{1-\frac{2}{d}}^{\frac{{\k}+1}{2}(d-2)} \\
&\ge &(d-1) \exp\left\{-\frac{2}{d-2}\frac{{\k}+1}{2}(d-2)\right\}\\
&=& (d-1) \exp\{-({\k}+1)\}\\
& > &1.
\end{eqnarray*}

\begin{lemma}\label{lem:variance}
Suppose that ${\k}\ge 1$ and $d$ is a sufficiently large positive even integer. Then,
$$
\frac{\E(H^2)}{\E(H)^2} \le (1+o(1)) \sqrt{\frac{d}{d-2({\k}+1)}}.
$$
\end{lemma}
Now Theorem \ref{thcolor} easily follows from this, since
$$\Pr(H=0)\leq  \frac{\Var(H)}{\E(H)^2}\leq  (1+o(1))\sqrt{\frac{d}{d-2({\k}+1)}}-1.$$

\subsubsection{Expectation (the proof of Lemma~\ref{lem:expectation})}

Let a $2m$-cycle in $\cX$ be a set of $2m$ disjoint pairs of points of $\cX$ such
that they form a $2m$-cycle
in $X$ (i.e. a Hamilton cycle) when they are projected by $\psi_1$ to $X$. Let
$p_{2m}$ be the probability that a given set of
$2m$ disjoint pairs of points of~$\cX$ forming a $2m$-cycle is contained in a random
configuration and that $\cU$ holds.

First note that from the proof of Lemma~\ref{spoil} the number of configurations partioned
into $2m$ cells of $d$ points for which $\cU$ holds is asymptotically
\begin{equation}\label{eq:no_loops}
\sim e^{-(d-1)/2}(2dm-1)!!=e^{-(d-1)/2}\frac{(2dm)!}{2^{dm}(dm)!}
\end{equation}
After fixing the pairs in a $2m$-cycle we have to randomly pair up $2(d-2)m$ points.
In other words, we want to compute the number of configurations partioned
into $2m$ cells of $(d-2)$ points for which $\cU$ holds. Hence, again by Lemma~\ref{spoil} we get,
\begin{equation*}
\sim e^{-(d-3)/2}(2(d-2)m-1)!!
\end{equation*}
and
$$p_{2m} \sim \frac{e^{-(d-3)/2}(2(d-2)m-1)!!}{e^{-(d-1)/2}(2dm-1)!!} =
e\frac{(2dm-4m-1)!!}{(2dm-1)!!}.
$$
Next, let $a_{2m}$ be the number of possible $2m$-cycles on~$\cX$. From (9.2) in~\cite{JLR} we get,
$$
a_{2m} = \frac{(d(d-1))^{2m} (2m)!}{4m}.
$$
Let $q_{2m}$ be the probability that a randomly chosen set $U$ of $2\k m$
points of $\cY$ (represented by $2m$ $\k$-sets) is equal (after the projection $\psi_2$)
to $\bX$, i.e., $\psi_2(U)=\bX$. Note that $U$ must contain precisely one copy
of every element of $\bX$. Hence, we have $(d/2)^{2\k m}$ out of $\binom{{\k}dm}{2{\k}m}$
choices for $U$. Thus, again by the proof of Lemma~\ref{spoil} we get,
\begin{equation*}
q_{2m} \sim \frac{e^{-(\k-1)(d-4)/4} (d/2)^{2\k m}}{e^{-(\k-1)(d-2)/4} \binom{{\k}dm}{2{\k}m}}
= e^{(\k-1)/2} \frac{(d/2)^{2\k m}}{\binom{{\k}dm}{2{\k}m}}.
\end{equation*}
Consequently,
\begin{align*}
\E(H) &= a_{2m} p_{2m} q_{2m} \\
&\sim e^{(\k+1)/2} \frac{d^{({\k}+1)2m}(d-1)^{2m} (2m)!
(2dm-4m-1)!!(2{\k}m)!({\k}dm-2{\k}m)!}{2^{2{\k}m+2} m (2dm-1)!! ({\k}dm)!}.
\end{align*}
Using the Stirling formula yields Lemma~\ref{lem:expectation}. Recall that
$(2N-1)!!\sim \sqrt{2}\bfrac{2N}{e}^N$.

\subsubsection{Variance (the proof of Lemma~\ref{lem:variance})}

Let $C_1$ and $C_2$ be two $2m$-cycles in $\cX$ sharing precisely $b$ pairs. Clearly, 
$|C_1\cup C_2|=4m-b$. Denote by $p_{2m}(b)$ the probability that $C_1$ and $C_2$ are 
contained in a random configuration of $\cX$ for which $\cU$ holds. (Clearly, $p_{2m}(2m)=p_{2m}$).
First note that if we ignore $\cU$ then the number of configurations containing $C_1$ 
and $C_2$ equals
$$
(2dm-2(4m-b)-1)!!
$$
Next conditioning on $\cU$ we obtain that the number of configurations containing $C_1$ 
and $C_2$ is bounded from above by
$$
e^{-(d-5)/2}(2dm-2(4m-b)-1)!!
$$
(The factor $e^{-(d-5)/2}$ corresponds to the case when $b=0$.)
Hence,
\begin{equation}\label{eq:no_loops_y}
p_{2m}(b) \le (1+o(1))\frac{e^{-(d-5)/2}(2dm-2(4m-b)-1)!!}{e^{-(d-1)/2}(2dm-1)!!} \sim
e^2\frac{(2dm-8m+2b-1)!!}{(2dm-1)!!}.
\end{equation}
Let $U$ and $W$ be two randomly chosen collections of $2m$ $\k$-sets in $\cY$ satisfying
$|W|=|U|=2m$ and $|W\setminus U|= 2m-b$.
Let $r_{2m}(b)$ be the probability
that both $U$ and $W$ are both equal (after the projection $\psi_2$) to $\bX$, 
i.e., $\psi_2(U)=\psi_2(W)=\bX$.
Conditioning on $\psi_2(U)=\bX$ we have $(d/2-1)^{2{\k}m-{\k}b}$ out of
$\binom{{\k}dm-2{\k}m}{2{\k}m-{\k}b}$ choices for $W$.
Thus, similarly as in \eqref{eq:no_loops_y} we obtain
$$
r_{2m}(b) \le (1+o(1))q_{2m} \frac{e^{-(\k-1)(d-6)/4} (d/2-1)^{2{\k}m-{\k}b}}{ e^{-(\k-1)(d-4)/4}
\binom{{\k}dm-2{\k}m}{2{\k}m-{\k}b}}
\sim e^{(\k-1)/2}q_{2m} \frac{ (d/2-1)^{2{\k}m-{\k}b}}{ \binom{{\k}dm-2{\k}m}{2{\k}m-{\k}b}}.
$$
Moreover, let $N(b)$ be the number of $2m$-cycles in $\cX$ that intersect a given $2m$-cycle in $b$
pairs.
By \cite{JLR} (cf. last equation on page~253), we get
$$
N(b) = \sum_{a=0}^{\min\{b,2m-b\}} \frac{2am}{b(2m-b)}2^{a-1}(d-2)^{2m+a-b}(d-3)^{2m-a-b}(2m-b-1)!
{b\choose a}{2m-b\choose a},
$$
where for $a=b=0$ we set $\frac{a}{b}=1$.

Consequently,
\begin{align*}
\frac{\E(H^2)}{\E(H)^2} &\le \frac{1}{\E(H)} +
\sum_{b=0}^{2m-1} \frac{N(b)p_{2m}(b)r_{2m}(b)}{a_{2m}p_{2m}^2q_{2m}^2} \\
&\le \frac{1}{\E(H)} +
(1+o(1)) \sum_{b=0}^{2m-1} \sum_{a=0}^{\min\{b,2m-b\}} \Bigg{(}
\frac{a(2m)^2}{b(2m-b)^2}2^a(d(d-1))^{-2m}(d-2)^{2m+a-b}\\
&\phantom{\le }\times (d-3)^{2m-a-b}  \binom{b}{a}\binom{2m-b}{a} 
\frac{(2m-b)! (2dm-8m+2b-1)!! (2dm-1)!!}{(2m)! (2dm-4m-1)!!^2}\\
&\phantom{\le }\times \frac{(d/2-1)^{2{\k}m-{\k}b}}{
\binom{{\k}dm-2{\k}m}{2{\k}m-{\k}b}} \frac{\binom{{\k}dm}{2{\k}m}}{(d/2)^{2{\k}m}}\Bigg{)}.
\end{align*}
Below we ignore all cases for which $a=0$, $a=b$ or $a+b=2m$
since their contribution is negligible as can be easily checked by the reader.
Using the Stirling formula, the terms in the sum can be written as
\begin{multline*}
\frac{1}{4\pi m} h(a/(2m),b/(2m))\exp\{2m\cdot g(a/(2m),b/(2m))\} \\
\times\left( 1+ O\left(\frac{1}{\min\{a,b-a,2m-a-b\}+1} \right)\right),
\end{multline*}
where
\begin{align*}
g(x,y) &= x\log(2) - \log(d) - \log(d - 1) + (1 + x - y)\log(d - 2) \\
&\qquad+ (1 - x - y) \log(d - 3) + y\log(y) + 2 (1 - y)\log(1 - y) \\
&\qquad- (y - x)\log(y - x) -  2 x\log(x) - (1 - x - y)\log(1 - x - y) \\
&\qquad+ (d/2 - 2 + y)\log(d - 4 + 2 y) + (d/2)\log(d) - (d - 2)\log(d - 2) \\
&\qquad+ {\k}(d/2 - 1)\log(d) + {\k}(1 - y)\log(1 - y) + {\k}(d/2 - 2 + y)\log(d - 4 + 2 y) \\
&\qquad- {\k}(d - 3 + y)\log(d - 2)
\end{align*}
and
$$
h(x,y) = \frac{ \sqrt{d(-4 + d + 2 y)} }{ {\sqrt{(d-2)^2 y(1-y)(1 - x - y)(y-x)}} }.
$$
Although the next computations may be verified by hand, the reader
might find the assistance of Mathematica useful. We give the definitions of $g(x,y)$ and $h(x,y)$ in
Mathematica format in Appendix~\ref{appendix}.

Now we analyze function $g(x,y)$ in the domain
$$
S=\{(x,y) : 0< x< y< 1-x\}.
$$
First, we compute the
first derivatives:
\begin{align*}
\frac{\partial g}{\partial x} &= \log(2)-\log(d-3) + \log(d-2) - 2 \log(x) + \log(-x + y)
+ \log(1 - x - y)\\
\frac{\partial g}{\partial y} &= -\log(d-3) - (1 + {\k})\log(d-2) - (2 + {\k})\log(1 - y)\\
&\qquad+ \log(1 - x - y) + \log(y) - \log(-x + y) + (1 + {\k}) \log(d-4 + 2 y).
\end{align*}
Let  $(x_0,y_0) = (2(d-2)/(d(d-1)), 2/d)$. Note that since
$\frac{\partial g}{\partial x}(x_0,y_0)=\frac{\partial g}{\partial y}(x_0,y_0)=0$, $(x_0,y_0)$
is a critical point of~$g$ and $g(x_0,y_0)=0$. Let $D^2{g}$ be the
Hessian matrix of second derivatives. Routine calculations show that
$$
D^2{g}(x,y)=\left(
\begin{array}{c c}
-\frac{2}{x} + \frac{1}{x - y} + \frac{1}{-1 + x + y} & \frac{1}{-x + y} + \frac{1}{-1 + x + y}\\
\frac{1}{-x + y} + \frac{1}{-1 + x + y}           &
\frac{2 + {\k}}{1 - y} + \frac{1}{x - y} + \frac{1}{y} + 
\frac{1}{-1 + x + y} + \frac{2 (1 + {\k})}{-4
+ d + 2 y}\\
\end{array}
\right)
$$
Hence,
$$
D^2{g}(x_0,y_0)=\left(
\begin{array}{c c}
-\frac{(d-1)^2 d}{2 (d-3)} & \frac{(d-4) (d-1)^2 d}{2 (d-2)(d-3)} \\
\frac{(d-4) (d-1)^2 d}{2 (d-2)(d-3)} & -\frac{d (16 + d (-34 + d (28 + (-9 + d) d - 2 {\k}) + 6
{\k})}{2 (d-3) (d-2)^2}
\end{array}
\right)
$$
One can verify that
$$
\Det(D^2{g}(x_0,y_0)) = \frac{d^3 (d-1)^2 (d - 2 (1 + {\k}))}{4 (d-3) (d-2)^2}.
$$
Since $-\frac{(d-1)^2 d}{2 (d-3)}<0$ and $\Det(D^2{g}(x_0,y_0))>0$ for $d>2(1+{\k})$,
we conclude that $D^2{g}(x_0,y_0)$
is negative definite at $(x_0,y_0)$. Hence, $g$ has a local maximum there.
Now we show that $(x_0,y_0)$ is the unique global maximum point of $g$ in 
$S$. Moreover, we argue that that $g(x,y)$ has no asymptote near the boundary 
of $S$, nor does it approach a limit which is greater than~$0$ (for $d$ large enough).

First recall that the function
\begin{equation}\label{eq:zlogz}
f(z) =
\begin{cases}
z\log(z) & \mbox{if } 0<z<1,\\
0 & \mbox{if } z=0 \mbox{ or } z=1
\end{cases}
\end{equation}
is continuous on $[0,1]$.
Consequently, function $g(x,y)$ can be extended to a continuous function on
$$
T=\{(x,y) : 0\le x\le y\le 1-x\}.
$$
Note that $-1/e \le f(z) \le 0$ (cf. \eqref{eq:zlogz}). Thus,
\begin{align*}
g(x,y) &\le \log(2) - \log(d) - \log(d - 1) + (1 + x - y)\log(d - 2) \\
&\qquad+ (1 - x - y) \log(d - 2) + 0 + 0 \\
&\qquad+ 1/e +  2/e + 1/e \\
&\qquad+ (d/2 - 2 + y)\log(d - 2) + (d/2)\log(d) - (d - 2)\log(d - 2) \\
&\qquad+ {\k}(d/2 - 1)\log(d) + 0 + {\k}(d/2 - 2 + y)\log(d - 2) \\
&\qquad- {\k}(d - 3 + y)\log(d - 2)\\
&= -y\log(d-2) + o(\log(d)),
\end{align*}
where the last term $o(\log(d))$ does not depend on $x$ and $y$.
Hence, there is a large enough $d$ such that $g(x,y)<0$ for all points in the
domain
$\{(x,y)\in T : 1/2(3+2{\k})\le y\}$.

Denote by $\partial T$ the boundary of $T$, i.e., $\partial T=T\setminus S$. In order to
finish, it is enough to show that:
\begin{enumerate}[(i)]
\item the only critical point in $\{(x,y)\in T\setminus\partial T : 
y\le 1/2(3+2{\k})\}$ is $(x_0,y_0)$, and\label{cond_1}
\item $g(x,y)<0$ for all points in $\{(x,y)\in \partial T : y\le 1/2(3+2{\k})\}$.\label{cond_2}
\end{enumerate}

Solving the equation $\frac{\partial g}{\partial y}(x,y)=0$ for $x$, noting that the
equation is linear in~$x$, we obtain
\begin{align*}
x = \frac
{y(1-y)\left( (d-3)(d-2)^{{\k}+1}(1-y)^{{\k}+1} -(d-4+2y)^{{\k}+1} \right)}
{(1-y)^{{\k}+2}(d-3)(d-2)^{{\k}+1} - y(d-4+2y)^{{\k}+1}}.
\end{align*}
Substituting this expression for~$x$ in $\frac{\partial g}{\partial x}(x,y)=0$
(actually in $\exp\{\frac{\partial g}{\partial x}(x,y)\}=1$) yields the equation
\begin{multline*}
0 = \psi(y) = 2(1-2y)^2 (1-y)^{\k} (d-4+2y)^{{\k}+1}(d-2)^{{\k}+2}
-y(1-y)^{2{\k}+2} (6-5d+d^2)^2 (d-2)^{2{\k}}\\
+2y(1-y)^{{\k}+1} (d-4+2y)^{1+{\k}} (d-3) (d-2)^{{\k}+1}
-y(d-4+2y)^{2+2{\k}}.
\end{multline*}
We see from our previous considerations that $\psi(y_0)=0$. It remains to show that
for large $d$, $y_0$ is the
only value in $\{y:0< y\le 1/2(3+2{\k})\}$ for which $\psi(y)=0$. To this
end we show that $\psi^\prime(y)<0$
implying that $\psi(y)$ is a monotone function (and clearly also continuous).
From the definition of $\psi(y)$ we get,
\begin{align*}
\psi^\prime(y) &=  \left(-y(1-y)^{2{\k}+2} (6-5d+d^2)^2 (d-2)^{2{\k}}\right)^\prime 
+ O(d^{2{\k}+3})\\
&= (1-y)^{2{\k}+1} (-1+y(2{\k}+3)) d^{2{\k}+4} + O(d^{2{\k}+3}),
\end{align*}
where the hidden constant in $O(d^{2{\k}+3})$ does not depend on $y$.
Hence, for a sufficiently large $d$ the derivative $\psi^\prime(y)<0$ for all 
$0< y\le 1/2(3+2{\k})$ (independently from $d$).
This shows that \eqref{cond_1} holds.

We split \eqref{cond_2} into three cases. One is for $0=x<y$, one for $0<x=y$ 
and the last one for $x=y=0$. Note that
\begin{align*}
g_1(y) &= g(0,y)\\ 
&= - \log(d) - \log(d - 1) + (1 - y)\log(d - 2) \\
&\qquad+ (1 - y) \log(d - 3) +  2 (1 - y)\log(1 - y) - (1 - y)\log(1 - y) \\
&\qquad+ (d/2 - 2 + y)\log(d - 4 + 2 y) + (d/2)\log(d) - (d - 2)\log(d - 2) \\
&\qquad+ {\k}(d/2 - 1)\log(d) + {\k}(1 - y)\log(1 - y) + {\k}(d/2 - 2 + y)\log(d - 4 + 2 y) \\
&\qquad- {\k}(d - 3 + y)\log(d - 2).
\end{align*}

Recall that $0<y<1/2(3+2{\k})$. It is easy to check that
$$
g_1^\prime(y) = -\log(d) + o(\log(d)),
$$
where the last term $o(\log(d))$ does not dependent on $y$.

Thus, for $d$ sufficiently large $g_1(y)$ is a decreasing function. Hence, by continuity
$$
g_1(y) \le g_1(0) = g(0,0).
$$

Later we show that $g(0,0)<0$.

Now let $0<x=y\le 1/2(3+2{\k})$. Define
\begin{align*}
g_2(y)&=g(y,y) \\
&= y\log(2) - \log(d) - \log(d - 1) + \log(d - 2) \\
&\qquad+ (1 - 2y) \log(d - 3) + y\log(y) + 2 (1 - y)\log(1 - y) \\
&\qquad-  2 y\log(y) - (1 - 2y)\log(1 - 2y) \\
&\qquad+ (d/2 - 2 + y)\log(d - 4 + 2 y) + (d/2)\log(d) - (d - 2)\log(d - 2) \\
&\qquad+ {\k}(d/2 - 1)\log(d) + {\k}(1 - y)\log(1 - y) + {\k}(d/2 - 2 + y)\log(d - 4 + 2 y) \\
&\qquad- {\k}(d - 3 + y)\log(d - 2).
\end{align*}
Consequently,
\begin{multline*}
g_2^\prime(y) = \log(2) - 2 \log(d-3) - \k \log(d-2) +  2 \log(1 - 2 y)\\ 
- (2 + \k) \log(1 - y) - \log(y) + (1 + \k) \log(d-4 + 2 y)
\end{multline*}
and
$$
g_2^{\prime\prime}(y) = (2 + \k)/(1 - y) - 1/y + 4/(-1 + 2 y) + 2 (1 + \k)/(d-4 + 2 y).
$$
Note that since $0<y\le 1/2(3+2{\k})$ we get that for $d$ large enough 
$g_2^{\prime\prime}(y)<0$. Thus, $g_2^\prime(y)$ is a decreasing function.
Moreover, since
$$
\lim_{y\to 0^{+}} g_2^\prime(y) = \infty
$$
and
$$
g_2^\prime(2/d) = 2 \log((d-4)/(d-3)) < 0,
$$
we conclude that $g_2(y)$ has a local maximum at $\x \in (0,2/d]$. Clearly 
such local maximum is the global maximum in the interval $(0,1/2(3+2{\k})]$. 
Unfortunately, it is not clear how to determine $\x$ since the equation 
$g_2^\prime(y)=0$ seems not to have any ``nice'' solution. Therefore, we define 
a new auxiliary function
$$
g_3(y) = g_2(y) - (2/3)  (d/2)^2  \log((d - 4)/(d - 3)) y^3
$$
on $(0,2/d]$. Clearly $g_2(y) \le g_3(y)$. Thus in order to show that $g_2(\x)<0$, 
it suffices to prove that $g_3(y)<0$ for any $y\in (0,2/d]$.
Analogously to analyzing $g_2(y)$ one can show that $g_3^{\prime\prime}(y)<0$ for $d$ 
large enough. Moreover, since $g_3^\prime(2/d)=0$, we get that $g_3(y)$ is an 
increasing function on $(0,2/d]$. Thus,
\begin{equation}\label{eq:g3}
g_3(y) \le g_3(2/d) = (8/3d - 1) \log((d-4)/(d-3)) + \log((d-2)/(d-1)).
\end{equation}
As one can check the right hand side in \eqref{eq:g3} is negative for sufficiently 
large $d$, as required.

It remains to show that $g(0,0)<0$. By continuity we get
$$
g(0,0) = \lim_{y\to 0^{+}} g_2(y) \le  g_3(2/d)<0.
$$

This completes the proof of \eqref{cond_2} and so the proof of showing that $(x_0,y_0)$
is the unique global maximum in~$T$.

The rest of argument is totally standard for such variance
calculations (see, e.g., \cite{FJMRW,JLR}). Finally, we obtain
\begin{align*}
\frac{\E(H^2)}{\E(H)^2} &\le (1+o(1)) \frac{1}{2\pi} \int_{-\infty}^{\infty}
\int_{-\infty}^{\infty} h(x_0,y_0)
\exp\left\{ -\frac{1}{2}(z_1,z_2)D^2{g(x_0,y_0)}(z_1,z_2)^T \right\} \ dz_1\ dz_2\\
&\sim \frac{h(x_0,y_0)}{\Det(D^2{g}(x_0,y_0))^{1/2}}\\
&= \frac{(d-1) d^2}{2(d-2)\sqrt{d-3}} \cdot \frac{2 (d-2) \sqrt{d-3}} { (d-1)
\sqrt{d^3 (d - 2 (1 + {\k}))}}\\
&=  \sqrt{\frac{d}{d-2({\k}+1)}},
\end{align*}
as required.

\section{Concluding remarks}
In this paper, we showed that $(\log n)/ n^{k-1}$ is the asymptotic threshold
for the existence of loose Hamilton
cycles in $H_{n,p;k}$ for $n$ a multiple of $2(k-1)$. It would be
nice to drop this divisibility requirement
and replace it by the necessary $(k-1) | n$, as mentioned in Introduction.
We address this question in our future work.

\section{Acknowledgment}
We would like to thank the anonymous referee for carefully reading this
manuscript, many helpful comments and pointing out some errors in the 
previous version of this manuscript. We are also very grateful to
Svante Janson and Nick Wormald for fruitful discussions about contiguity
of random regular graphs (contiguity was used in the previous version of this paper).

\appendix
\section{Mathematica expressions}\label{appendix}
For convenience, we replace here $\k$ by $k$.
{\small{
\begin{verbatim}
g[x_,y_,d_,k_] = x Log[2] - Log[d] - Log[d - 1] + (1 + x - y) Log[d - 2] \
     + (1 - x - y) Log[d - 3] + y Log[y] + 2 (1 - y) Log[1 - y] \
     - (y - x) Log[y - x] -  2 x Log[x] - (1 - x - y) Log[1 - x - y] \
     + (d/2 - 2 + y) Log[d - 4 + 2 y] + (d/2) Log[d] - (d - 2) Log[d - 2] \
     + k(d/2 - 1) Log[d] + k(1 - y) Log[1 - y] + k(d/2 - 2 + y) Log[d - 4 + 2 y] \
     - k(d - 3 + y) Log[d - 2];
\end{verbatim}
}}
{\small{
\begin{verbatim}
h[x_,y_,d_] =  Sqrt[d(-4 + d + 2 y)] / Sqrt[(d-2)^2 y(1-y)(1 - x - y)(y-x)];
\end{verbatim}
}}

\end{document}